\documentclass[12pt]{amsart}
\usepackage{times,amsfonts,amsmath,amstext,amsbsy,amssymb,
amsopn,amsthm,upref,eucal
}
\newcommand \Ker {{\rm Ker}}
\newcommand \R {{ \mathbb R}}
\newcommand \RR {{ \mathbb R}}
\newcommand \si {{ \sigma}}

\newcommand \T {{ \mathbb T}}
\newcommand \E {{ \mathcal E}}
\newcommand \cL {{ \mathcal L}}
\newcommand \g {{ \frak g}}
\newcommand \B {{ \mathcal B}}
\newcommand \Z {{ \mathcal Z}}
\newcommand \Vi {{ \rm Vect^\infty(M)}}
\newcommand \betal {{ \Psi}}

\newcommand \GN {{ \Gamma\setminus N}}

\newcommand \CGN {{ C^\infty(\Gamma\setminus N)}}

\newcommand \reps{{
(\Gamma\setminus N\hat )}}

\newcommand \ci{{C^\infty}}

\def \eps{ \varepsilon}
\def \al{ \alpha}
\def \la{\lambda}
\def \th{\theta}

\def \om{\omega}

\newtheorem {thm} {Theorem}[section]
\newtheorem {prop} {Proposition}[section]
\newtheorem{cor}{Corollary}[section]

\newtheorem{definition}{Definition}[section]
\newtheorem{conj}{Conjecture}[section]

\newenvironment{pf}{ Proof:}{ $\square$ \newline}

\begin{document}

\title[Actions with globally hypoelliptic leafwise Laplacian]{Actions with globally hypoelliptic leafwise Laplacian and rigidity}
\author[Danijela Damjanovi\'c]{Danijela Damjanovi\'c }
\thanks{ Based on research supported by NSF grant DMS-1001884  and NSF grant DMS-1150210}

\subjclass[2013]{}
\address{Department of Mathematics, Rice University, 6100 Main Str. Houston, TX, 77005}

\email{dani@rice.edu}  
%\{http://www.math.psu.edu/damjanov}

\begin{abstract} In this paper we prove several results concerning smooth $\mathbb R^k$ actions with the property that their leafwise Laplacian is globally hypoelliptic. Such actions are necessarily uniquely ergodic and minimal, and cohomology  is often finite-dimensional, even trivial.  Further we consider a class of examples of $\mathbb R^2$ actions on 2-step nilmanifolds, which have globally hypoelliptic leafwise Laplacian,  and we show transversal local rigidity under certain Diophantine conditions. 
\end{abstract}

\maketitle

\section{Introduction} 

\subsection{Transversal local rigidity of actions}
Let $\rho$ be  a smooth $\R^k$ action by diffeomorphisms  of a smooth compact manifold $M$. We say that $\rho$ is \emph{transversally locally rigid} if there exists a finite dimensional transversally $C^1$ family $\{\rho_\la\}_{\la\in \mathbb R^d}$ of smooth $\R^k$ actions on $M$, such that $\rho_0=\rho$, and such that every sufficiently small perturbation of the family $\rho_\la$ in a neighborhood of $\la=0$, intersects the smooth conjugacy class of $\rho$, or of an action obtained from $\rho$ via a small coordinate change in the acting group. In such a situation we also say that $\rho$ is \emph{transversally locally rigid with respect to the family $\{\rho_\la\}$}. This kind of local behavior for higher rank actions was studied in \cite{DK_UNI} for certain homogenenous parabolic actions and was further discussed in \cite{D}.

%(!!!) NOt what I want! Action $\rho$ is \emph{parameter locally rigid} if there exists a finite dimensional transversally $C^1$ family $\{\rho_\la\}_{\la\in \mathbb R^d}$ of smooth $\R^k$ actions on $M$ such that $\rho_0=\rho$ and such that every sufficiently small perturbation of $\rho$ is in the smooth conjugacy class of some  $\rho_\la$ for $\la$ close to 0.(!!!) Remove sec 6.4??????? Or discuss why the result should not be true for ANY perturbation.

\subsection{Background: Diophantine vector fields on tori, rigidity and cohomology} 

Let $M=\T^{n}$. Let  $\{\rho_\la\}$ be the family of constant vector fields: for $\la=(\la_1, \dots, \la_n)\in \R^n$, $\rho_\la=\sum_{k=1}^n \la_n\frac{\partial}{\partial x_k}$. We call the vector field $\rho_\la$  Diophantine if $\la$ is a Diophantine vector:  $|\la\cdot m|\ge C\|m\|^{-\gamma}$, for some constants $C>0$, $\gamma>n-1$, and for all $m\in \mathbb Z^n\setminus \{0\}$.

Recall the well known 
%\foot{check in Herman or Zehnder the details}
 result on perturbations of Diophantine vector fields: 

\begin{thm}[V.I. Arnold, J. Moser]\label{KAMtorus}
For $\om$ Diophantine, there is $l>0$ and $\epsilon>0$ such that for every smooth perturbation $\tilde\rho_\om=\rho_\om+\beta$ of $\rho_\om$ which is $\eps$-close to $\rho_\om$ in $C^l$-norm, there exists $\bar\la$ and there exists a smooth diffeomorphism $h$ of $\T^n$ such that:
$$dh^{-1}(\rho_{\bar\la}+\beta)\circ h=\rho_\om$$
\end{thm}

One can formulate this theorem in terms of families of perturbations as follows: 

{\it Let $\om$ be Diophantine. Then there exist $l>0$ and $\eps>0$ such that for every $f$ which is $\eps$-small in  $C^l$-norm, the family $\{ \rho_\la+\beta\}_{\la\in \R^n}$ of perturbations of $\rho_\om$, intersects the conjugacy class of $\rho_\om$ at $\rho_{\bar\la}+\beta$ for some parameter $\bar\la$. }

In fact, the same conclusion holds for \emph{any} transversal family of perturbations: 

\begin{thm}[\cite{D}]\label{KAMtorus-fam}
For $\om$ Diophantine, $\rho_\om$ is transversally locally rigid with respect to the family $\{\rho_\la\}$. Namely, there is $l>0$ and $\epsilon>0$, and a neighborhood $B$ of $\om$ in $\R^n$,  such that for every family of smooth perturbations $\{\tilde\rho_\la\}=\{\rho_\om+\tilde \beta_\la\}$ of  $\{\rho_\om\}$ which is $\eps$-close to $\{\rho_\la\}$ in $C^l$-norm for all $\la\in B$, and in $C^1$ norm transversally,  there exists $\bar\la\in B$ and there exists a smooth diffeomorphism $h$ of $\T^n$ such that $\tilde\rho_{\bar\la}$ is in the conjugacy class of $\rho_\om$:
$$dh^{-1}(\tilde\rho_{\bar\la})\circ h=\rho_\om$$
\end{thm}

Theorem \ref{KAMtorus-fam} has a simple proof which uses Theorem \ref{KAMtorus}, but it can also be derived from the main result of \cite{D} (See the first section of \cite{D} for more details). Both Theorem \ref{KAMtorus} and Theorem \ref{KAMtorus-fam} are about $\mathbb R$-actions, i.e. flows. For actions of larger groups it turns out that the result in Theorem \ref{KAMtorus-fam} is perhaps more natural, and certainly is more general. Namely, perturbing a family of actions with a single small perturbation like in Theorem \ref{KAMtorus} is not possible any more because it would result in a family which would not necessarily consist of abelian actions.

Consider a smooth $\R^k$, $k\ge 2$, action $\rho$ on a torus $M=\T^n$. It is defined by a $k$-tuple of commuting vector fields. One example is obtained by taking a $k$-tuple of constant vector fields $\la_1, \dots, \la_k$ (defined in standard coordinates by the k-tuple of vectors in $\R^n$, which we also denote by $\la_i$). 
Let $\{\rho_\la\}$ be the family of all $\mathbb R^k$ smooth actions on $\mathbb T^n$ which are generated by  such $k$-tuples of constant vector fields. Let  $\th=(\th_1, \dots, \th_k)$ be simultaneously Diophantine, i.e. $\max_i\{\|m\cdot \th_i\|\}\ge C\|m\|^{-\gamma}$ for some $C>0$, $\gamma>n-1$ and for all $m\in \mathbb Z^{n}$. Then the following is obtained in \cite{D}: 
%as a corollary of the general result from \cite{D}:

\begin{thm}[\cite{D}]\label{commuting_torus}
The action $\rho_\theta$ for $\theta$ simultaneously Diophantine is transversally locally rigid with respect to the family $\{\rho_\la\}$. 
\end{thm}
% in a neighborhood of a simultaneously Diophantine $\la=(\la_1, \dots, \la_k)$ is transversaly locally rigid at any Diophantine  $\th=\la$. This fact is also a corollary of the main result  in \cite{D2}.  

All the above local results rely crucially on the fact that the first cohomology over the actions involved is trivial, i.e, it reduces to constant cocycles. The space of cocycles over a smooth action can be viewed as the formal tangent space to the space of smooth actions in the neighborhood of the given action, thus its trivialization suggests strongly that the local picture about the given action should be simple as well. 

Also, the actions considered in the above theorems have another property (which in the case of a flow coincides with having trivial cohomology), namely, they have globally hypoelliptic leafwise Laplacian. It turns out that this property for a general, not necessarily homogeneous action, is very strong, and among other facts it also implies \emph{finite dimensional} first cohomology. Again, this strongly suggests that for such actions one may hope to obtain some classification of perturbations. 

\subsection{Results}
 In the first part of this paper  we recall what global hypoellipticity is,  and we prove several general results concerning dynamical and cohomological properties of $\mathbb R^k$ actions with globally hypoelliptic leafwise Laplacian: unique ergodicity, minimality and that first and second cohomology tend to be finite dimensional, in some cases even trivial. 

In the second part of the paper we generalize the local results described in Theorems \ref{KAMtorus-fam} and \ref{commuting_torus}  to 2-step nilmanifolds. Namely we prove transversal local rigidity for certain $\mathbb R^2$ homogeneous actions on 2-step nilmanifolds which satisfy certain Diophantine conditions. This is Theorem \ref{main} in Section \ref{main-section}. The proof of Theorem \ref{main} uses a general result from \cite{D} that describes cohomological conditions which imply transversal local rigidity (Theorem \ref{IFT}). The fact that these cohomological conditions are satisfied for the actions on 2-step nilmanifolds considered in Theorem \ref{main} comes from global hypoellipticity of the leafwise Laplacian and from the analysis of the induced action in nilmanifold representations. 

In Section \ref{further_remarks} we discuss some further questions related to actions with globally hypoelliptic leafwise Laplacian. 

Author would like to express gratitude to Giovanni Forni and Federico Rodriguez-Hertz for useful discussion and insights.

\section{Preliminaries}
\subsection{Globally hypoelliptic differential operators} Let $L$ be a differential operator on a compact connected smooth manifold $M$.  Denote by $\E$ the space $\ci(M)$ of smooth functions and by $\E'$ the space of distributions. Let $\Ker L=\{D\in \mathcal E': D(Lf)=0 \,\, \mbox{for all}\,\, f\in \E\}$. In what follows we assume that $\Ker L$ contains a nowhere vanishing smooth volume form on $M$. 

\begin{definition} Operator $L$ is globally hypoelliptic (GH) if every distributional solution $u$  to the equation $Lu=f$ where $f\in \E$, is in $\E$.  
\end{definition}

In \cite{GW1} Greenfield and Wallach prove the following fact about GH operators:
\begin{prop}[\cite{GW1}]\label{prop1}
If differential operator $L$ is GH, then $\dim (\Ker L)<\infty$.
\end{prop}
Greenfield and Wallach  apply the above fact to $L=\cL_X$, the Lie derivative for a given smooth vector field $X$ on $M$, and they obtain:
\begin{prop}[\cite{GW2}]\label{prop2} 
If $L=\cL_X$ is GH, then $\dim (\Ker L)=1$ and $\Ker L$ is spanned by a nowhere vanishing smooth volume form on $M$.
\end{prop}
This immediately implies the following fact for the flow $\Phi_X^t$ tangent to a GH vector field $X$: 
\begin{cor}[\cite{GW1}]
If $L=\cL_X$ is GH, then $\Phi_X^t$ is uniquely ergodic and minimal flow on $M$.
\end{cor}
%This already significantly restricts possibilities for a flow tangent to a  GH vector field $X$. 

Prime example of  a flow tangent to a GH vector field is the  Diophantine flow  on a torus. An open question is

\begin{conj}\label{GWK}[Greenfield-Wallach, A. Katok]
If $L=\cL_X$ is GH then $\Phi_X^t$ is smoothly conjugated to the Diophantine flow on the torus. 
\end{conj} 
Reader interested in recent developments in the direction of this conjecture may consult \cite{F}. The formulation above is due to Greenfield and Wallach. A. Katok's conjecture involved a seemingly stronger property of $\Phi_X^t$ being cohomology free.
Operator $L=\cL_X$ (or the corresponding flow $\Phi_X^t$) is \emph{cohomology free} if the equation $Lu=f$ has a $\ci$ solution for every $\ci$ function $f$ whose average on $M$, with respect to the nowhere vanishing smooth volume form, is zero. Equivalence of the two notions in case of the Lie derivative operator for a smooth vector field $X$ was proved in \cite{CC}, and a complete proof can also be found in \cite{F}:

\begin{prop}[\cite{CC}, \cite{F}]
$L=\cL_X$ is GH if and only if it is cohomology free. 
\end{prop}

The proof in \cite{F} of the above Proposition actually shows that for any linear operator $L$ which is GH and for which $\Ker L$ is trivial i.e. $\Ker L$  consists of constant functions, the range $\rm Im(L)$ of $L:C^\infty (M)\to C^\infty (M)$, is a closed space. This implies by the Hahn-Banach theorem that the smooth cohomology over $L$ is completely determined by the space of $L$-invariant distributions. If this space is 
trivial, then $L$ is cohomology free. Otherwise, if for a GH operator $L$, $\Ker L$ is not trivial and if the range $\rm Im(L)$ of $L$ is a closed space, then by Proposition \ref{prop1}, $\rm Im(L)$ must be in the intersection of \emph{finitely} many distributions 
%finite dimensional in case $L$ is GH, it is spanned by finitely many distributions 
$D_1, \dots, D_m$. Thus we have: 

\begin{prop}\label{cohL}
a) If $L$ is GH and $\rm Im(L)$ is a closed space, then there are finitely many distributions $D_1, \dots, D_m$ such that 
%$\ker L=<D_1, \dots, D_m>$, then the following holds: 
for every smooth $f\in \cap_{i=1}^m \Ker D_i$,  there exists a function $h$ in $\ci(M)$ such that $Lh=f$. 

b) If $L$ is GH and $\Ker L$ consists of constant functions, then $\rm Im(L)$ is a closed space and  $L$ is cohomology free. 
\end{prop}

{\it Remark.} In \cite{GW1} it is also proved that if $L$ is GH and commutes with an elliptic operator, then it does have a closed range. 

%\begin{cor}
%If $L$ is GH and $\Ker L$ is trivial (i.e. contains only constant functions), then $L$ is cohomology free.  
%\end{cor}

\subsection{Cohomology over $\R^k$ actions}
%Let  $\g$ be a finite dimensional Lie algebra and $V$ a $\g$-module with representation $X\mapsto \L_{\rho(X)}$, where $\rho:\g\to V$ and $ \L_{\rho(X)}$ is an endomorphism of $V$ induced by  $\rho(X)\in V$, $X\in \g$.  Let $\frak C^k(\frak g, V)= \rm Hom_k(\wedge^k\g, V)$, $k=0, 1, \dots, \rm dim (\frak g)$, be the space of cochains,  i.e. the vector space of multilinear alternating mappings from $\g\times\dots\times\g$ to $V$. In this context, $\frak C^0(\g, V)=V$.

Let $\rho$ be a smooth action of $\R^k$ on manifold $M$, generated by $k$ commuting smooth vector fields $X_1, \dots, X_k$ on $M$. Let $\g$ denote the linear space spanned by $X_1, \dots, X_k$. Let $V$ denote either $\ci(M)$, or the space of smooth vector fields $\Vi$, or in the case when $M$ is paralellizable, the space of vector fields with constant coefficients.  

Denote by  $\frak C^j(\frak g, V)$ the space $ \rm Hom_j(\wedge^j\g, V)$  of multilinear alternating mappings from $\g\times\dots\times\g$ to $V$. This is the space of j-cochains where  $j=0, 1, \dots, k$,  Clearly $\frak C^0(\g, V)=V$. The coboundary operators $\delta_\rho^{j}: \frak C^j(\frak g, V)\to \frak C^{j+1}(\frak g, V)$ are defined by:
%where for $\omega\in \frak C^k(\frak g, V)$
\begin{equation*}
\begin{aligned}
(\delta_\rho^{j}\omega)(Y_0, Y_1, \dots, Y_j)&=\sum _{0\le i\le j}(-1)^j\cL_{Y_i}\omega(Y_0, \dots, \hat Y_i, \dots, Y_j)\\
%&+\sum_{0\le i\le l\le j}(-1)^{i+l}w([[X_i, X_l], \dots, \hat X_i, \dots, \hat X_l, \dots, X_j)
\end{aligned}
\end{equation*}
for $\omega\in \frak C^j(\frak g, V)$, $Y_0, \dots, Y_j\in \g$, and $ \hat Y_i$ means that the $i$-th entry is omitted.  

Denote by $\Z^j_\rho(V)=\Ker (\delta_\rho^{j})$ the space of all $j$-cocycles and by $\B^j_\rho(V)=\rm Im(\delta_\rho^{j-1})$ denote the space of all $j$-coboundaries. The factor space $\Z^j_\rho(V)/\B^{j}_\rho(V)$ is the $j$-th cohomology over the action $\rho$ with coefficients in $V$, and is denoted by $H_\rho^j(V)$. 
%This  is the usual Chevalley-Eilenberg cohomology of a Lie algebra, with coefficients in the module $V$.

We will be interested only in the first two coboundary operators, and here is how they look like:
\begin{equation*}
\begin{aligned}
&(\delta_\rho^0 H) (X)=\cL_{X}H, \,\,\, H\in V, \,\,X\in \g\\
&(\delta_\rho^1 \betal)(X, Y)=\cL_{X}\betal(Y)-\cL_{Y}\betal(X),\betal\in \frak C^1(\frak g, V)
\end{aligned}
\end{equation*}

It is easy to see that for an $\mathbb R^k$ action $\rho$, $H^{j}_\rho(V)$ is trivial for $j>k$, and that $\Ker \delta_\rho^k$ is the whole space of $k$-cochains.

Left or right sided inverses of operators $\delta_\rho^j$, if they exist, will be denoted by $\delta_\rho^{j*}$. 

%(1) Let $N$ be a compact smooth manifold. Let $ V=\Vi N$ and let $\rho:\g\to \Vi $ be a Lie algebra homomorphism defining  a smooth $G$-action. Taking the Lie derivative $\L_{\rho(X)}$  with respect to $\rho(X)\in \Vi M$ for $X\in \g$, turns $V$ into a $\g$-module. Given a sequence of seminorms $\{\|\cdot\|_k\}$ (for example  $C^k$ norms) which makes $V$ a tame Fr\'echet space, notice that the space $\frak C^1$ is also a  tame Fr\'echet space with the sequence of seminorms defined as the maximum of $\|\cdot\|_k$ norms  of images of  a fixed basis of $\g$ under $\rho$. Similarily, any  $\frak C^k$ is a tame Fr\'echet space. The coboundary operators $\delta_\rho^0$ and $\delta_\rho^1$ in this case are $(\mathcal C, 1)$-tame, where the constants $C_k$ in the collection $\mathcal C=\{C_k\}$ depend on the $C^k$ norms of $\rho$. 

%(2) Let $V=\frak h$ be a Lie algebra, and $\rho:\g\to \h$ a Lie algebra homomorphism. Let $\L_{\rho(X)}$ be defined by the Lie bracket in $\h$: $\L_{\rho(X)}(H)=[\rho(X), H]$. This cohomology appears when we consider actions by left translations on homogeneous spaces, such as in \cite{DKpar} and \cite{KW}.  

%(3) Special case of the above is $V=\frak g$ and $\rho:\g\to \g$ is a Lie algebra automorphism. We will be interested in the case when $\rho$ is the identity map, i.e the cohomology $H^1_{\rm Id}(\g, \g)$.
%\subsection{Spaces, norms and tame operators}

\section{Properties of smooth volume preserving $\R^k$ actions with globally hypoelliptic leafwise Laplacian}
Assume an $\R^k$ action $\rho$ on $M$ is generated by $k$ smooth commuting vector fields $X_1, \dots X_k$. Let $$L_\rho=\sum_{i=1}^kX_k^2$$ be the leafwise Laplacian associated to the action $\rho$, where $X_k$ stands for the Lie derivative $\cL_{X_k}$. (Remark: It would certainly be more standard to call $-\sum_{i=1}^kX_k^2$ the leafwise Laplacian, but change of the sign has no influence on any statements in this paper, so we opted for a simpler form). 

\subsection{Unique ergodicity}

\begin{prop} If $L_\rho$ is GH, then $\dim(\cap_{i=1}^k \Ker X_i))=1$.
\end{prop}
\begin{pf} The argument here is similar to the proof of \ref{prop2} in \cite{GW2}.\\
Since clearly $\cap_{i=1}^k \Ker X_i\subset \Ker L_\rho$, from Proposition \ref{prop1} it follows that $\dim(\cap_{i=1}^k \Ker X_i)<\infty$. Surely the nowhere vanishing smooth form $\omega$ is in $\cap_{i=1}^k \Ker X_i$. Assume that for some $f\in \ci(M)$ we also have $f\omega\in\cap_{i=1}^k \Ker X_i$. This implies $X_if=0$ for $i=1, \dots, k$, and since $X_i$'s are vector fields we also have $X_if^n=0$ for any $n\ge 0$. But then $\cap_{i=1}^k \Ker X_i$ contains $\{\omega, f\omega, f^2\omega, \dots\}$. Since $\cap_{i=1}^k \Ker X_i$ is finite dimensional this means that for some $N$ large enough we must have $\sum_{j=0}^Na_jf^j=0$. But this implies that $f$ takes only finitely many values, and since $M$ is connected, it follows that $f$ is constant. 
\end{pf}

An immediate consequence is:
\begin{cor}
If $L_\rho$ is GH, then $\rho$ is uniquely ergodic and minimal $\mathbb R^k$ action. 
\end{cor}

\emph{Remark.} Flows or $\mathbb R^k$ actions whose distributional kernel contains only constant multiples of a smooth volume form are sometimes called \emph{strongly uniquely ergodic} actions.

\subsection{ Finite dimensional cohomology}
\begin{prop}\label{1coh-trivial}
a) If $L_\rho$ is GH and if $\rm Im (L_\rho)$ is a closed space, then $\dim H^1_\rho(\ci(M))< \infty$. 
%(How about $H_\rho^k(\ci)$ for general k???????)

b) If $L_\rho$ is GH and $\Ker L_\rho$ is 1-dimensional spanned by the invariant volume form, then $ H^1_\rho(\ci(M))$ is trivial, namely every cocycle with trivial average is a coboundary. 
\end{prop}
\begin{pf}
%CHECK THE WHOLE PROOF

\emph{Proof of part a).} It is clear that $\omega$ in $\in \frak C^1(\g, C^\infty)$ is completely determined by its values on the basis of $\g$, i.e. by $\omega(X_i)$, $i=1, \dots, k$. 

Define operator $d:  \frak C^1(\g, C^\infty)\to C^\infty _0(M)$ by 
$$d(\omega)=\sum_{i=1}^kX_i\omega(X_i).$$ 

Also define $L(\omega)$ to be the element of $\frak  C^1(\g, C^\infty)$ such that  $$(L(\omega))(X_i)=L_\rho(\omega(X_i)).$$ 

Then the  following properties hold:

(1) $d\circ\delta^0_\rho=L_\rho$ on $C^\infty_0(M)$.

(2) $L=\delta^0_\rho\circ d$ on $\Ker \delta^1_\rho$.

The property (1) is trivial to check. For property (2), given $\omega$, let $f_i:=\omega(X_i)$. Then we  have: $$\delta^0_\rho\circ d(f_1, \dots, f_k)=\delta^0_\rho(\sum_i X_if_i)=(\sum_i X_i X_1f_i, \dots, \sum_i X_i X_k f_i)$$ Now use the assumption that $\omega\in \Ker \delta^1_\rho$, i.e.  that $X_jf_i=X_if_j$. By using this substitution in the previous formula we get: $$\delta^0_\rho\circ d(f_1, \dots, f_k)=(\sum _iX_i^2f_1, \dots, \sum_i X_i^2 f_k)=L(\omega)$$ This proves property (2).

Now let $\omega\in \frak C^1(\g, C^\infty)$ satisfy $\delta^1_\rho(\omega)=0$, and consider the smooth function $d(\omega)$. 

Since $L_\rho$ is GH, by Proposition \ref{prop1} it follows that $\Ker L_\rho$ is finite dimensional and smooth, say it is spanned by $F_0, F_1, \dots, F_m$, where $F_0=1$.  Denote by $D_i$ the distributions induced by $F_i$. 

Assume first that $d(\omega)\in \cap_{i=0}^m \Ker D_i$. Since $\rm Im (L_\rho)$ is a closed space by assumption, by Proposition \ref{cohL} a), there exists $h\in C^\infty (M)$ such that 

$$L_\rho h=d(\omega)$$
Because of property (1), this implies $d(\delta^0_\rho h)=d(\omega)$, therefore $d(\delta^0_\rho h-\omega)=0$.  Since $\delta^0_\rho h-\omega\in \Ker \delta^1_\rho$, by the property (2) we have  

$$L(\delta^1_\rho h-\omega)=\delta^1_\rho( d(\delta^1_\rho h-\omega))=0$$
Thus for each $i$: $L(\delta^0_\rho h-\omega)(X_i)=L_\rho((\delta^0_\rho h-\omega)(X_i))=0$

Since  $\Ker L$ is finite dimensional spanned by $F_0, \dots, F_m$, we have for each $i=1,\dots, k$ some constants $a_j^i$ such that:

$$(\delta^0_\rho h-\omega)(X_i)=\sum_{j=0}^m a^i_j(\omega)F_j.$$
So the (finite dimensional) cohomology of $L$ completely determines the cohmology classes of $\omega$ for which $d(\omega)$ is in the kernel of invariant distributions for $L_\rho$. If we denote by $\gamma_\omega$ the form belonging to a finite dimensional space of forms, defined by $ \gamma_\omega(X_i)=\sum_j a^i_j(\omega)F_j$, then for each $i=1,\dots, k$:
\begin{equation}\label{omega-coh}
\omega(X_i)=X_i h-\gamma_\omega(X_i)
\end{equation}

In order to treat the case when $d(\omega)$ is not in $\cap \Ker D_i$,  define functionals $$\tilde D_i:=D_i\circ d$$ on $\frak C^1(\g, C^\infty)$. The space $\Ker \delta^1_\rho$ is a closed subspace in $\frak C^1(\g, C^\infty)$, so for every non-trivial $\tilde D_i$ there exists a form $\alpha_i$ in $\Ker \delta^1_\rho$ such that $\tilde D_i(\alpha_j)=\delta_{ij}$. 

Now if $d(\omega)\notin \cap \Ker D_i$, then for constants  $b_i(\omega):=\tilde D_i(\omega)$, clearly $$d(\omega-\sum_{i=0}^m b_i(\omega)\alpha_i)\in \cap_{i=0}^m \Ker D_i $$  
so the previous part of the proof applies, giving that there exists a smooth $h$ such that
$$\omega-\sum_{i=0}^m b_i(\omega)\alpha_i=\delta^0_\rho h-\gamma_\omega$$
This gives a decomposition of $\omega$ into a coboundary and a form whcih belongs to a fixed finite dimensional space of forms:
\begin{equation}\label{omega-coh2}
\omega=\delta^0_\rho h+ \sum_{i=0}^m b_i(\omega)\alpha_i -\sum_{j=0}^m a^i_j(\omega)F_j
\end{equation}
Therefore, modulo a cocycle belonging to a finite dimensional subspace of  $\frak C^1(\g, C^\infty)$, $\omega$ belongs to the range of the operator $\delta^0_\rho$.  This implies that the first cohomology is finite dimensional. The dimension of the cohomology $H^1_\rho(\ci(M))$ is at most $2m+1$, where $m+1$ is the dimension of the first cohomology  for $L_\rho$.

\emph{Proof of part b).} From Proposition \ref{cohL} b) the range $\rm Im (L_\rho)$ is closed, thus the proof of part a) applies only with a significant simplification in the beginning: assume that $\Ker (L_\rho)$ is spanned by $F_0=1$. Then $D_0(d(\omega))=0$ for all $\omega$ so only the first part of the proof of a) is relevant, moreover in \eqref{omega-coh} $\gamma_\omega$ takes constant values, so modulo a constant cocycle $\omega$ is a coboundary, i.e.  $H^1_\rho(C^\infty(M))$ is trivial.

 \end{pf}
 
 \emph{Remark.} There certainly exist actions with trivial $H^1_\rho(C^\infty(M))$ which do not have GH leafwise Laplacian. One example is the Weyl chamber flow on $SL(n, \mathbb R)/\Gamma$ for $n\ge 3$. (See for example \cite{KS}).

\begin{prop}\label{2coh-trivial}
If $L_\rho$ is (GH) and $\rm Im (L_\rho)$ is closed, then the top cohomology $H_\rho^k(\ci(M))$ is finite dimensional and  has the dimension less or equal than the kernel of $L_\rho$. In particular, $H_\rho^k(\ci(M))$ is trivial if $L_\rho$ has kernel consisting only of constants. 
\end{prop}

\begin{pf} For the top cohomology, 
any k-cochain $\omega$ is completely defined by the value it takes on  $X_1, \dots, X_k$.  Also, in the top cohomology $\Ker \delta_\rho^k$ is the whole space of $k$-cochains, i.e. smooth functions. \\
Let  $s=\omega(X_1, X_2, \dots, X_k)\in C^\infty (M)$. Let as before $D_0, D_1, \dots, D_m$ be the distributions spanning $\Ker L_\rho$. If $s\in \cap_ {i=0}^m \Ker D_i$, then by Proposition \ref{cohL} there exists $h\in \ci$ such that $$s=Lh=X^2_1h+\dots+X^2_kh=X_1(X_1h)+\dots+X_k(X_kh).$$ Thus if we define a $k-1$ chain $\phi$ by: $$\phi(X_1, \dots, \hat X_i, \dots, X_k):=(-1)^{i+1}X_ih,$$ then clearly $\delta_\rho^k\phi=\omega$. For every one of the distributions $D_i$ there exists $s_i$ such that $D_i(s_j)=\delta_{ij}$, so we have that for every $k$-cochain $\omega$: $$\omega(X_1, \dots, X_k)-\sum_{i=0}^m D_i(\omega(X_1, \dots, X_k))s_i$$ lies in $\cap_ {i=0}^m \Ker D_i$ and therefore is a $k$-coboundary. Thus the top cohomology is finite dimensional of the dimension at most $m+1$. 
\end{pf}

\section{Transversal local rigidity  for $\R^2$ actions on 2-step nilmanifolds}\label{main-section}

\subsection{Setup and the formulation of the result} 
Let $\frak N$ be a 2-step rational nilpotent Lie algebra and let $N$ be the corresponding connected simply connected Lie group. Let $\Gamma$ be a (cocompact) discrete subgroup of $N$ and let $M=\Gamma\setminus N$.

Denote by $Y_1\dots, Y_q$ and $Z_1, \dots, Z_p$ a linear basis for  $\frak n$ (selected from $\log \Gamma$) so that:

$\bullet$ $Y_1 +[\frak N, \frak N], \dots, Y_q +[\frak N, \frak N]$ is a basis for $\frak N/[\frak N, \frak N]$, and

$\bullet$ $Z_1, \dots, Z_p$ is a basis for $[\frak N, \frak N]$.

Now define a $p+q$-dimensional family of $\mathbb R^2$ actions $\rho_{a, b}$: for  vectors $a=(a_1, \dots, a_q)\in \mathbb R^q$ and $b=(b_1\dots, b_p) \in   \mathbb R^p$,  action $\rho_{a, b}$ is generated by the commuting pair 
$$X_1:=a_1Y_1+\dots +a_qY_q,\, \,  \makebox{and}\,\, X_2:=b_1Z_1+\dots +b_p Z_p.$$

Recall that a vector $a=(a_1, \dots, a_n)\in \mathbb R^n$ is Diophantine, if there exist positive constants $C$ and $\gamma>n-1$ such that 
$$|a_1k_1+\dots +a_nk_n|\ge C\|(k_1, \dots, k_n)\|^{-\gamma}$$ 
holds for any non-trivial integer vector $ (k_1, \dots, k_n)$.

We consider perturbations of the $p+q$-dimensional family of $\mathbb R^2$ actions  $\{\rho_{a, b}\}$ in a neighborhood of $\rho_{\al, \beta}$ where $\al$ and $\beta$ are both Diophantine. 

The following is the perturbation result which will be proved in the subsequent sections:

\begin{thm}\label{main}  If $\al$ and $\beta$ are Diophantine vectors,  then $\rho_{\al,\beta}$ is transversally locally rigid with respect to the family $\{\rho_{a,b}\}_{a\in \R^q, b\in R^p}$.

%The family $\{\rho_{a,b}\}_{a\in \R^q, b\in R^p}$ on the two step nilmanifold is transversally locally rigid in neighborhood of every element $\rho_{\al,\beta}$ of the family, for which $\al$ and $\beta$ are Diophantine vectors. 
\end{thm} 

The starting point towards the proof the Theorem \ref{main} is the following result proved by Cygan and Richardson: 

\begin{thm}[Cygan-Richardson \cite{CR}]\label{CR} If $\al$ and $\beta$ are Diophantine, then for the action $\rho_{\al, \beta}$ the leafwise Laplacian $L_ {\rho_{\al, \beta}}$ is GH. Also, $\Ker L_{\rho_{\al, \beta}}$ is trivial, i.e. consists only of constant functions. 
%(DOUBLECHECK: correct from CR, but may need extra explanation??)
\end{thm}

\begin{cor}\label{trivial-coh}
a) If $\al$ and $\beta$ are Diophantine, then the action $\rho_{\al, \beta}$ has trivial first and second cohomology with coefficients in smooth functions.

b)   If $\al$ and $\beta$ are Diophantine, then the first cohomology over the action $\rho_{\al, \beta}$ with coefficients in vector fields coincides with the constant cohomology over the action, namely the first cohomology with coefficients in $\frak N$. 
\end{cor}
\begin{pf}
Part a) is follows by  Propositions \ref{1coh-trivial}  and \ref{2coh-trivial}. Part b) follows from part a) and Proposition 1 from \cite{D}. Namely, given that $\frak N$ is nilpotent, the first and second coboundary operators on vector fields have an upper triangular form (with respect to a basis chosen in $\frak N$), so after a finite number of inductive steps a) implies b). 
\end{pf}

%\begin{pf}The proof reduces to checking the conditions of the Theorem \ref{IFT}. Conditions (i) and (ii) are proved in Sections \ref{1coh-inverse}, Section \ref{2coh-inverse} and Section \ref{vector-fields}. It remains to check the condition (iii).  Part a) of condition (iii) is proved in Section \ref{constant-coh}.
%\end{pf}

The proof of Theorem \ref{main} is an application of the general "implicit function" type theorem proved in \cite{D}. This general result lists cohomological conditions on a finite dimensional  family of actions in a neighborhood of a given action in the family,  which imply transversal local rigidity. We state this general result  in Section \ref{general-result}. Two main cohomological conditions for family $\{\rho_{a,b}\}$ are already partially confirmed by Corollary \ref{trivial-coh}, however in addition to cohomology trivialization we will also need to obtain \emph{tame} estimates for the inverses $\delta_{\rho_{\al, \beta}}^{0*}$ and $\delta_{\rho_{\al, \beta}}^{1*}$, of the first two coboundary operators over the action $\rho_{\al, \beta}$. 

The $j-th$ coboundary operator has a  \emph{tame} inverse, if  for some scale of norms $\|\cdot\|_r$ (in this paper we use Sobolev norms), the following is satisfied for a smooth $j$-cochain $\Psi$: $$\|\delta_{\rho_{\al, \beta}}^{j*}\Psi\|_r\le C_r\|\Psi\|_{r+\sigma},$$ where $r>0$ is arbitrary,  $C_r$ is a constant which besides $\rho_{\al, \beta}$ and $\GN$ may only depend on $r$, and $\sigma$ is a constant which depends only on $\rho_{\al, \beta}$ and $\GN$. 

The tame estimate for the inverse of the first coboundary operator is obtained in Section \ref{1coh-inverse}, and the tame estimate for the inverse of the second coboundary operator is obtained in Section \ref{2coh-inverse}. For the purpose of obtaining tame estimates we need to make use of elementary Kirillov theory describing representations of nilpotent manifolds, and we summarize the facts from Kirillov theory needed here in Section \ref{kirilov}.  Finally, we show where the family $\rho_{a,b}$ comes from: namely that it is precisely the family which, together with coordinate changes for the action $\rho_{\al,\beta}$, parametrizes the constant cohomology i.e. the first cohomology over $\rho_{\al,\beta}$ with coefficients in $\frak N$. This computation is contained in the Section \ref{constant-coh}. The proof of Theorem \ref{main} is contained in Section \ref{proof-main}.

%For the purpose of the main result of this paper, we prove a stronger fact: the first and second cohomology over Diophantine $\rho_{\al, \beta}$  is in fact trivial (i.e. reduces to constants). Moreover,  the first and second  coboundary operators have tame inverses. We extend these facts to cohomology with coefficients in vector fields. This is the content of the next 5 sections. We use the same approach as in \cite{C-R} and \cite{FF}, that is we use Kirilov theory,  except that here more precision is needed in order to  obtain \emph{tame} estimates for solutions of coboundary equations. For readers' convenience in the next section we outline the basics of Kirilov theory, especially its consequences when restricted to the case of 2-step nilmanifolds considered in this paper. 

\subsection{Representations and Kirilov theory}\label{kirilov}

In the subsequent sections, given a differential operator $L$ the main goal is to study regularity of solutions of equation $Lu=f$ where $f\in C^\infty(\Gamma\setminus N)$. Classical approach to this type of a problem is to use irreducible Fourier series decomposition of a distribution $u$. Existence and properties of this decomposition for functions of nilmanifolds come from Kirilov theory \cite{K}, \cite{CGbook}. 

%Kirilov theory provides a correspondence between the irreducible unitary representation of $N$ and the coadjoint orbits in the space $\frak n^*$ of functionals on $\frak n$.  These are orbits of the co-adjoint action...

When $\frak N$ is  a finite dimensional real nilpotent Lie algebra and $N=\exp \frak N$ the corresponding Lie group, it was proved by Malcev that $N$ has a cocompact discrete subgroup $\Gamma$ if and only if $\frak N$ has rational structure constants with respect to some suitable basis. In this case the basis for $\frak N$ can be selected from $\log \Gamma$. 

Denote by $\hat N$ the space of equivalence classes of irreducible representations of $N$. It is proved by Kirillov \cite{K}, that the elements in $\hat N$ are in one-to-one correspondence with the coadjoint orbits in the linear dual $\frak N^*$ of $\frak N$. These are orbits under the coadjoint action on $\frak N$ which is defined by $Ad^*(g)\la=\la\circ Ad(g^{-1})$ for $g\in N$ and $\la\in \frak N^*$. Representation $\pi$ which corresponds to the coadjoint orbit $\mathcal  O(\la)=Ad^*(N)\la$ of $\la\in \frak N^*$ can be denoted by $\pi_\la$ and the orbit  $\mathcal  O(\la)$ may also be denoted sometimes by $\mathcal  O(\pi)$ if it corresponds to $\pi$. 

Denote the subspace of $\hat N$ which appears  in the discrete direct sum decomposition of $L^2(\Gamma\setminus N)$ by $\reps$. Richardson proved that $\pi\in \reps$ if and only if there exists $\la\in O(\pi)$ and a rational subalgebra $\frak M$ subordinate to $\la$ such that $\chi_\la(\exp \frak M\cap \Gamma)=1$. This implies that if $\frak Z$ is the center of $\frak N$, $\la$ must be integer valued on $\frak Z\cap \log \Gamma$. 

Given a representation $\pi$ corresponding to the coadjoint orbit $\mathcal O(\pi)=\mathcal O(\la)$, $\|\pi\|$ denotes the distance from $\mathcal O(\pi)$ to the origin in $\frak N^*$. 

%Denote by $\mathcal H_\pi$ the $\pi$-primary summand in $L^2(\Gamma\setminus N)$. 

Every nilmanifold $\GN$ has the associated torus $\Gamma[N, N]\setminus N$. The only representations in $\reps$ which are not infinite dimensional are the one-dimensional representations of $\Gamma[N, N]\setminus N$. For $f\in C^\infty (\GN)$, $f_0$ denotes the sum of these one dimensional components of $f$, and the set of the remaining infinite dimensional representations will be denoted by $\reps_\infty$. On a 2-step nilpotent Lie group $N$ the representations in $\reps_\infty$, rather the corresponding coadjoint orbits, can be parametrized by hyperplanes in the real vector space of the dimension equal to the dimension of the basis of $\frak N$. 

Theorem \ref{CR} in \cite{CR} uses Schr\"odinger models for the representations of two-step nilmanifolds in  $\reps_\infty$. Here we only need the fact that for every coadjoint orbit $O(\la)$ corresponding to an irreducible unitary representation in $\reps_\infty$,  one can assign a (unitarily equivalent) representation $\pi$ for which $d\pi$ acts on the center $\frak Z$ non-trivially by scalar multiplication: 
$$d\pi(Z)=i\la(Z)I$$
for every $Z\in [\frak N, \frak N]$,  \cite{CGbook}. Also, for a two-step nilmanifold $\GN$ the representations space $\mathcal H_\pi$ corresponding to any $\pi\in \reps_\infty$ is $L^2(\mathbb R^n)$, where $n$ varies with the structure and the dimension of $\frak N$ . 
%to every coadjoint orbit corresponding to a representation in $\reps_\infty$, one can assign a (unitarily equivalent) representation $\pi$ for which $d\pi(Z)$ acts as a multiplication by $i\la(Z)$ for every $Z\in [\frak N, \frak N]$, \cite{GWbook-ionaj drugipapir}.

%{\it DO THIS PART BETTER, USE GEOMETRIC DESCRIPTION FROM CHENG !!!!!!!!!!!!!!!!!!!!!!!!!!!!!!!!!!!!: There are two types of irreducible representations for a 2-step $N$: one kind of representations are those which are trivial on the center and they are coming from the torus $\Gamma[N, N]\setminus N$; the other kind are representations which are not trivial on the center. These are infinite dimensional representations and by Stone von Neumann theorem they can be described by Schr\"odinger representations. } 

For non-toral (infinite dimensional) representations $\pi$, from the geometric description of the corresponding coadjoint orbit,  in \cite{CHENG} the distance  $\|\pi\|$ was computed in terms of the values of the corresponding functional $\la$  on the basis of $\frak N$. Immediate consequence is that: 
$$\|(\la(Z_1), \dots, \la(Z_p))\|\le \|\pi\|$$
where $Z_1, \dots, Z_p$ is the basis of the center of the two-step nilpotent Lie algebra $\frak N$.

%----------don't need the part bellow--------
%Then $\mathcal H_\pi$ has an irreducible decomposition $\mathcal H_\pi=\mathcal H_{\pi, 1}+\dots+ \mathcal H_{\pi, m(\pi)}$ where $m(\pi)$ is the (finite) multiplicity of $\pi$. 

%Function $f\in L^2(\Gamma\setminus N)$ decomposes as $$f=\sum_{\pi\in \reps} f_\pi$$ and $f_\pi$ has further decomposition as $$f_\pi=\sum_{i=1}^{m(\pi)} f_{\pi, q}$$ 

%It was proved by Ausleander and Brezin that $f\in C^\infty (\Gamma\setminus N)$ implies that each $ f_{\pi, q}$ is $C^\infty$ and that both sums $f=\sum f_\pi=\sum f_{\pi, q}$ converge uniformly for a $C^\infty$ function $f$. Conversely, for $C^\infty$ functions $ f_{\pi, q}$, $f=\sum f_{\pi, q}$ is $C^\infty $  if and only if $\sum _{\pi, q}\|Uf_{\pi,q}\|_2<\infty$ for each $U$ in the universal enveloping algebra of $\frak N$.  

%---------------------

Every $f\in L^2(\GN)$ decomposes as: 

$$f=\sum _{\pi\in \reps}f_\pi$$
where $f_\pi$ denotes the component  of $f$ in the representation $\pi$. It is well known that if $f\in \CGN$ then $f_\pi$ are smooth. Also Sobolev space $W^r( \GN)$ decomposes into a direct sum of Sobolev spaces $W^r(\mathcal H_\pi)$ where $\pi\in \reps$ and $\mathcal H_\pi$ is the representation space of $\pi$.

We will use here a more precise statement which is obtained by Corwin and Greenleaf in \cite{CG}, and which compares the Sobolev norms of components $f_\pi$ to the Sobolev norm of $f$, for a smooth function $f$. Namely the following holds: 

\begin{thm}[\cite{CG}]\label{CG}
a) For a smooth $f$ for any $s>0$ and $k>0$:  $$\|f_\pi\|_s\le C_s\|\pi\|^{-k}\|f\|_{s+k}$$ where $C_s$ is a constant depending only on $s$ and the manifold. 

b) There exists exponent $k_0>0$ such that for $k\ge k_0$: $$\sum_{\pi\in \reps} \|\pi\|^{-k}<\infty.$$
\end{thm}

%Besides obtaining smooth solutions to cohomological equations, we will also need to obtain estimates for the solutions, and this will be done in each irreducible representation showing up in $\reps$. For this purpose we discuss now Sobolev spaces....

\emph{Remark.} In the subsequent sections we will use letter $C$ to denote \emph{any} constant which depends only on the manifold $\GN$ and the action $\rho_{\al, \beta}$, and $C_r$ do denote additional dependence of the constant on parameter $r$.

\subsection{Tame inverse   $\delta_\rho^{0*}$ to the first coboundary operator}\label{1coh-inverse}

%Q: WHY DOES THIS PF NOT WORK ON GENERAL NILMANIFOLDS?

Recall that $\delta_\rho^{0}$ is the first coboundary operator which to every smooth function $h$ assigns a cocycle on the Lie algebra $\frak g$ of the acting group, defined by $X\mapsto \mathcal L_{X} h$ for $X\in \frak g$. This cocycle is clearly completely determined by its values on the basis of the acting Lie algebra, and for the action $\rho=\rho_{\alpha, \beta}$ these are precisely the vector fields $X_1$ and $X_2$ which generate the action. 

Thus,  $\delta_\rho^{0}$ can be considered as a map $h\mapsto (\mathcal L_{X_1}h, \mathcal L_{X_2} h)$. For more compact notation denote $\mathcal L_{X_i}$ by $X_i$, $i=1,2$. From Corollary \ref{trivial-coh}, we know that  $\delta_\rho^{0}$ has an inverse operator $\delta_\rho^{0^*}$ on cocycles of trivial average. However, for the application of Theorem \ref{IFT}, one also needs \emph{tame} estimates for $\delta_\rho^{0^*}$. The goal of this section is to obtain such estimates.  

Denote by $C_0^\infty(\GN)$ the space of smooth functions on $\Gamma\setminus N$ with zero average. Let $f, g\in C_0^\infty(\Gamma\setminus N)$ generate a smooth cocycle over $\rho_{\al, \beta}$. Since from Corollary \ref{trivial-coh} we already have that the first cohomology over $\rho_{\al, \beta}$ trivializes, it follows that there exists $h\in \CGN$ such that:

\begin{equation}\label{1coh}
\begin{aligned}
X_1 h&=f\\
X_2h&=g
\end{aligned}
\end{equation}
We will obtain estimates for $h$ by looking at  \eqref{1coh} in toral and non-toral representations. 

If $\pi$ denotes an irreducible representation in $\reps$,  then as before $f_\pi$ denotes the component of $f$ in the representation space of $\pi$. Every function $f$ decomposes as: 
$$f=f_0+ \sum_{\pi\in \reps_{\infty}} f_\pi$$
where $f_0$ lives on the associated torus $\Gamma[N, N]\setminus N$, and  $\reps_{\infty}$ denotes the set of non-toral (infinite dimensional) representations in $\reps$. 

%we denote by $f_\pi$, $g_\pi$ and $h_\pi$ the components of $f, g$ and $h$ in this representation, respectively.  

Let $\pi$ be an irreducible representation associated to the coadjoint orbit of $\la\in \frak N^*$, and let $d\pi$ denote the corresponding representation of the Lie algebra  $\frak N$. Recall that $\la([\frak N, \frak N]\cap \log \Gamma)\subset \mathbb Z$. 

The coboundary equation \eqref{1coh} in representation $\pi$ looks like: 
\begin{equation}\label{1coh-rep}
\begin{aligned}
d\pi(X_1) h_\pi&=f_\pi\\
d\pi(X_2)h_\pi&=g_\pi
\end{aligned}
\end{equation}

For the toral component $f_0$, the first equation above is:
$$\sum_i\al_i\frac{\partial}{\partial x_i}h_0=f_0$$
This is the well known small divisor equation on the associated torus $\Gamma[N, N]\setminus N$, and it is a classical fact  (for example \cite{Moser}) that given the Diophantine condition on $\al$, the norm of $h_0$ is estimated by: 
\begin{equation}\label{est-toral}
\|h_0\|_r\le C_r\|f_0\|_{r+\sigma}
\end{equation}
for any positive $r$, where $C_r$ and $\sigma$ are fixed constants depending only on Diophantine properties of $\al$ and the dimension of $\Gamma[N, N]\setminus N$.

Note that for the toral component $f_0$ of a smooth function $f$ we have: 
\begin{equation}\label{tor-all-est}
\|f_0\|_r\le C_r\|f\|_{r+\sigma'}
\end{equation}
 where $\sigma'$ depends on the dimension of $\GN$ and the dimension of the associated torus. This is an immediate consequence of Theorem \ref{CG} a), and the fact that that toral representations are parametrized by integer vectors. 

Now if $\pi$ belongs to the class of non-toral representations, then the second equation in \eqref{1coh-rep} becomes:

$$i\sum_j\beta_j\la(Z_j)h_\pi=g_\pi$$
Using the Diophantine condition on $\beta$, this implies:
\begin{equation}\label{est-nontoral}
\|h_\pi\|_r\le C\|(\la(Z_1), \dots, \la(Z_p))\|^\gamma\|g_\pi\|_r
\end{equation}
for any positive $r$, where $C$ and $\gamma$ are positive constants which depend on the Diophantine properties of $\beta$.

Since $\|(\la(Z_1), \dots, \la(Z_p))\|\le \|\pi\|$,  it follows that for every $\pi$ the estimate for $h_\pi$ of the following form holds: 

\begin{equation*}\label{1coh-est}
\|h_\pi\|_r\le C\|\pi\|^\gamma \|g_\pi\|_{r}
\end{equation*}
where $r$ is arbitrary positive,  and constants $C$ and $\sigma$ depend on the Diophantine properties of $\al$ and $\beta$. 

From the Theorem \ref{CG} of Corwin-Greenleaf:
 $$\|g_\pi\|_{r}\le C_r\|\pi\|^{-k}\|g\|_{r+k}$$
  Putting together the estimates above it follows that $h$ satisfies:
  
\begin{equation}\label{1coh-est-final}
\|h\|_r\le C_r\|f\|_{r+\sigma+\sigma'}+ C_r\sum_{\pi\in \reps_\infty} \|\pi\|^{-k+\gamma}\|g\|_{r+k}
\end{equation}
From Theorem \ref{CG} part b),  there exists $k_0$ such that $\sum_{\pi\in \reps} \|\pi\|^{-k+\gamma}<\infty$ as long as $k-\gamma\ge k_0$, therefore, by choosing $k$ to be the nearest integer greater than $\gamma +k_0$, and by redefining $\sigma:=\max\{\sigma+\sigma', k\}$ we obtain a tame estimate for $h$:
\begin{equation}\label{1coh-est}
\|h\|_r\le C_r\max\{\|f\|_{r+\sigma}, \|g\|_{r+\sigma}\}
\end{equation}
where $C$ and $\sigma$ are constants depending only on the manifold $\GN$ and the Diophantine properties of $\al$ and $\beta$.

%\subsection{Tame solutions for the leafwise Laplacian equation}???Don't need any more??

\subsection{Tame inverse $\delta_\rho^{1*}$ to the second coboundary operator }\label{2coh-inverse} 

%IN THIS SECTION FIX THE CONSTANTS TO REFLECT CONSTANT COHOMOLOGY

Second coboundary operator $\delta_\rho^{1}$ takes a mapping from the acting Lie algebra to the smooth functions, into a mapping from the double product of the acting Lie algebra $\frak g\times \frak g$ to the smooth functions. Given that our acting Lie algebra is two dimensional  and generated by vector fields $X_1$ and $X_2$, the operator $\delta_\rho^{1}$ can be simply described as a map from $\CGN\times \CGN$ into $\CGN$, defined by  $(f, g)\mapsto X_2f-X_1g$. In this section we show the existence of a \emph{tame} inverse for this operator, defined on the range of the operator. 

Let $f, g, \phi\in \CGN$ be such that 
\begin{equation}\label{2coh}
X_2f-X_1g=\phi
\end{equation}

If $\pi$ is an irreducible representation (corresponding to $\la\in \frak N^*$), then the second cohomology equation \eqref{2coh} in this representation looks like:

\begin{equation}\label{2coh-reps}
d\pi(X_2)f_\pi-d\pi(X_1)g_\pi=\phi_\pi
\end{equation}
Since $d\pi(X_2)$ acts trivially, for the full toral components $f_0$, $g_0$ and $\phi_0$, equation \eqref{2coh-reps} reduces to: 
$$-\sum \al_i\frac{\partial}{\partial x_i}g_0=\phi_0$$
As before, this is a standard small divisor equation on the associated torus, so using the Diophantine assumption on $\alpha$, and \eqref{tor-all-est}, we obtain the  following estimate:
\begin{equation}\label{est-toral-g}
\|g_0-g_{triv}\|_r\le C_r\|\phi_0\|_{r+\sigma}\le C\|\phi\|_{r+\sigma+\sigma'}
\end{equation}
where $g_{triv}$ denotes the component of $g$ in the trivial representation, and $C$, $\sigma$ and $\sigma'$ are constants depending on Diophantine properties of $\alpha$. Also, as in the previous section, under the condition that $f_0$ is 0 in the trivial representation, we can solve the equation $f_0=d\pi(X_1)h_0=\sum_i \al_i\frac{\partial}{\partial x_i} h_0$ for $h_0$ with the estimate:
$$\|h_0\|_r\le C_r\|f\|_{r+\sigma+\sigma'}$$
So on the associated torus we can define $\tilde f_0:=0$,  $\tilde g_0:= g_0-g_{triv}$, $f_{triv}$ to be the component of $f$ in the trivial representation,  and this gives the following splitting of the toral components:
\begin{equation}
\begin{aligned}
f_0&=X_1h_0+\tilde f_0+f_{triv}\\
g_0&=X_2h_0+\tilde g_0+g_{triv}
\end{aligned}
\end{equation}
with estimates:
\begin{equation}\label{est-toral}
\begin{aligned}
\max\{\|\tilde f_0\|_r, \|\tilde g_0\|_r\}&\le C_r\|\phi\|_{r+\sigma+\sigma'}\\
\|h_0\|_r&\le C_r\max\{\|f\|_{r+\sigma+\sigma'}, \| g\|_{r+\sigma+\sigma'}\}
\end{aligned}
\end{equation}

If $\pi$ is a non-toral representation, then as in the previous section we can solve the equation $d\pi(X_2)\bar g_\pi=g_\pi$. Namely, $\bar g_\pi=-i(\sum_j\beta_j\la(Z_j))^{-1}g_\pi$, which, after using the Diophantine assumption on $\beta$, implies the following estimate for $\bar g_\pi$:
$$\|\bar g_\pi\|_r\le C\|(\la(Z_1), \dots, \la(Z_p))\|^{\gamma}\|g_\pi\|_{r}$$
Similarily, let $\bar \phi_\pi$ be the solution to the equation $d\pi(X_2)\bar \phi_\pi=\phi_\pi$, satisfying the estimate: 
$$\|\bar \phi_\pi\|_r\le C\|(\la(Z_1), \dots, \la(Z_p))\|^{\gamma}\|\phi_\pi\|_{r}$$
Then using commutativity of $d\pi(X_1)$ and $d\pi(X_2)$, from \eqref{2coh-reps} we have:
$$f_\pi=d\pi(X_1)\bar g_\pi +\bar \phi_\pi.$$

Thus by defining $\tilde f_\pi:=\bar\phi_\pi$ and $\tilde g_\pi:=0$, the $\pi$-components of $f$ and $g$ can be split as: 

\begin{equation}
\begin{aligned}
f_\pi&=d\pi(X_1)\bar g_\pi+\tilde f_\pi\\
g_\pi&=d\pi(X_2)\bar g_\pi+\tilde g_\pi
\end{aligned}
\end{equation}
with estimates:
\begin{equation}\label{est-non-toral}
\begin{aligned}
\max\{\|\tilde f_\pi\|_r, \|\tilde g_\pi\|_r\}&\le C\|(\la(Z_1, \dots, \la(Z_p))\|^{\gamma}\|\phi_\pi\|_{r}\\
\|\bar g_\pi\|_r&\le  C\|(\la(Z_1, \dots, \la(Z_p))\|^{\gamma}\max\{\|f_\pi\|_r, \| g_\pi\|_r\}
\end{aligned}
\end{equation}

We may now summarize information obtained from all representations.  As before, it will be also used that $\|(\la(Z_1, \dots, \la(Z_p))\|\le \|\pi\|$. Let $f, g$ and $\phi$ satisfy the equation \eqref{2coh}. Then we can split $f$ and $g$ into a part  which is a coboundary and into a part which is of the order of $\phi$ in the following way: 

\begin{equation}\label{splitting}
\begin{aligned}
f&=X_1H+\tilde f+f_{triv}\\
g&=X_2H+\tilde g+g_{triv}
\end{aligned}
\end{equation}
where

$$H:= h_0+\sum_{\pi\in \reps_\infty}\bar g_\pi$$

$$\tilde f:= \tilde f_0+ \sum_{\pi\in \reps_\infty}\tilde f_\pi$$

$$\tilde g:= \tilde g_0+ \sum_{\pi\in \reps_\infty}\tilde g_\pi$$
with $\tilde f_0=0$, $\tilde g_0=g_0-g_{triv}$, $\tilde f_\pi=\bar \phi_\pi$ and $\tilde g_\pi=0$, for $\pi\in \reps_\infty$.

From all the previously obtained estimates \ref{est-toral} and \ref{est-non-toral}, it follows that: 
\begin{equation*}
%\begin{aligned}
\|H\|_r\le  C_r\max\{\|f\|_{r+\sigma+\sigma'}, \| g\|_{r+\sigma+\sigma'}\}+ C_r\sum_{\pi\in \reps_\infty}\|\pi\|^{\gamma-k}\max\{\|f\|_{r+k}, \| g\|_{r+k}\}
%\end{aligned}
\end{equation*}
and 
\begin{equation*}
\max\{\|\tilde f\|_r, \|\tilde g\|_r\}\le C_r\|\phi\|_{r+\sigma+\sigma'}+C_r\sum_{\pi\in \reps_\infty}\|\pi\|^{\gamma-k}\|\phi\|_{r+k}
\end{equation*}

Again, by choosing a fixed $k>k_0+\gamma$, (where $k_0$ is the constant from Theorem \ref{CG} b)), and by redefining $\sigma:=\max\{\sigma+\sigma', k\}$, the final tame estimates follow:

\begin{equation}\label{splitting-est}
\begin{aligned}
\|H\|_r&\le  C_r\max\{\|f\|_{r+\sigma}, \| g\|_{r+\sigma}\}\\
\max\{\|\tilde f\|_r, \|\tilde g\|_r\}&\le C_r\|\phi\|_{r+\sigma}
\end{aligned}
\end{equation}

It is clear from the above splittings and estimates that if $(f, g)$ is a cocycle, i.e. if $\delta_\rho^1((f, g))=0$,  then $\phi=0$ and thus modulo the averages of $f$ and $g$, the cocycle $(f, g)$ is a coboundary.

{\bf Remark}. The above is a direct proof.  However due to the fact that the 2nd cohomology  over the action $\rho:=\rho_{\al, \beta}$ is also trivial,  there is an alternative way of obtaining the tame estimates for $\delta^{1*}_{\rho}$, which is perhaps a more canonical approach in the setting of actions with GH leafwise Laplacian. Namely, we have from Corollary \ref{cohL}  that the cohomology over $L_{\rho}$ is trivial, so  for any $C^\infty$ function $s$ if the average of $s$ is trivial,  then there exists a $C^\infty$ function $h$ such that $L_{\rho}h=s$. (Notice that since $\ker L_{\rho}$ contains only constants, the solution $h$ is unique up to a constant). Suppose that in addition one can obtain tame estimate for $h$ with respect to $s$. (In the context of action $\rho$ it is possible to prove this fact by looking into each irreducible representation: in toral representation the estimate is the classical small divisor problem, and in non-toral representation the estimate comes from solving a second order differential equation.) 
Now this  implies $s=X_1^2h+X_2^2h=X_2 \tilde  f-X_1 \tilde g$ where $\tilde f:=X_2h$ and $\tilde g=-X_1h$.
This would clearly  define a tame inverse $\delta_\rho^{1*}$ for the operator $\delta_\rho^{1}$ on the full space of smooth 2-cochains (which are in fact smooth functions). If $s=X_2f-X_1g$ then $f-\tilde f$ and $g-\tilde g$ constitute a coboundary, so there exists $h' $ such that $f-\tilde f=X_1 h'$ and $g-\tilde g=X_2 h'$. Thus in this version, the splitting for $f$ and $g$ would look  like:
$$f=X_1h'+X_2h$$
$$g=X_2h'-X_1h$$
where $h'$ satisfies tame estimates with respect to $(f, g)$ and $h$ satisfies tame estimates with respect to $s=X_2f-X_1g$. Notice that in this version of splitting, $h$ and $h'$ are unique up to a constant.

%\subsection {$\delta_\rho^{0*}$ and $\delta_\rho^{1*}$ for vector fields.}\label{vector-fields}
%Section 5.3, Proposition 1  in \cite{D2}.

\subsection{Computation of constant cohomology}\label{constant-coh}

Consider an arbitrary constant cocycle $\betal$ over  $\rho=\rho_{\al, \beta}$. This is a map from the abelian Lie algerba $\frak g$ spanned by $X_1$ and  $X_2$, into the nilpotent Lie algebra $\frak N$. It is completely defined by $\betal(X_1)$ and $\betal(X_2)$. Recall that $\dim \frak N=p+q=n$. Let 
$$\betal(X_k)=\sum_{i=1}^q a^k_iY_i +\sum _{j=1}^p b^k_j Z_j$$
where $k=1$ or $2$. 

If $\betal\in \Ker \delta_\rho^1$, then $[X_2,\betal(X_1)]=[X_1,\betal(X_2)$], that is:

$$\sum_{i=1}^q a^1_i[X_2, Y_i] +\sum _{j=1}^p b^1_j[X_2,  Z_j]=\sum_{i=1}^q a^2_i[X_1,Y_i] +\sum _{j=1}^p b^2_j [X_1, Z_j]$$

Recall that $X_1:=\sum_{l=1}^q\al_lY_l$ and $X_2:=\sum_{t=1}^p\beta_tZ_t$. Therefore:

$$\sum_{i=1}^q a^1_i\sum_{t=1}^p\beta_t[Z_t, Y_i] +\sum _{j=1}^p b^1_j\sum_{t=1}^p\beta_t[Z_t,  Z_j]=\sum_{i=1}^q a^2_i\sum_{l=1}^q\al_l[Y_l,Y_i] +\sum _{j=1}^p b^2_j \sum_{l=1}^q\al_l[Y_l, Z_j]$$
Using the bracket conditions in $\frak N$ the above equality becomes:
$$0=\sum_{i=1}^q \sum_{l=1}^q a^2_i\al_l[Y_l,Y_i] $$
which implies that
$$a_i^2\alpha_l-a_l^2\alpha_i=0$$ 
for all $i\ne l$ is the condition for $\betal$ being in $\Ker \delta_\rho^1$. This practically means that the vector $(a_1, \dots, a_q)$ is proportional to the vector $(\al_1,\dots, \al_q)$. 

Now we investigate the condition for $\betal\in \rm{Im} \delta_\rho^0$. Suppose $H= \sum_{i=1}^q h_iY_i +\sum _{j=1}^p h'_j Z_j$ and assume $\betal=\delta^1H$. This implies 

\begin{equation}
\begin{aligned}
\betal(X_1) =[X_1, H]&
=\sum_{i=1}^q h_i[X_1, Y_i] +\sum _{j=1}^p h'_j [X_1, Z_j]\\
&=\sum_{i=1}^q h_i\sum_{l=1}^q\al_l[Y_l, Y_i] +\sum _{j=1}^p h'_j \sum_{l=1}^q\al_l[Y_l, Z_j]\\
&=\sum_{i=1}^q h_i\sum_{l=1}^q\al_l[Y_l, Y_i] =\sum_{i=1}^q \sum_{l<i} (h_i\al_l-\al_ih_l)[Y_l, Y_i]
\end{aligned}
\end{equation}
where $[Y_l, Y_i]\in \{Z_1, \dots, Z_p\}$.
Similarily
\begin{equation}
\begin{aligned}
\betal(X_2) =[X_2, H]&
=\sum_{i=1}^q h_i[X_2, Y_i] +\sum _{j=1}^p h'_j [X_2, Z_j]\\
&=\sum_{i=1}^q h_i\sum_{t=1}^p\beta_t[Z_t, Y_i] +\sum _{j=1}^p h'_j \sum_{t=1}^p\beta_t[Z_t, Z_j]=0\\
\end{aligned}
\end{equation}
This implies that the condition for $\betal$ (defined by $\betal(X_k)=\sum_{i=1}^q a^k_iY_i +\sum _{j=1}^p b^k_j Z_j$, $j=1,2$) to be in $\rm{Im}\delta_\rho^0$ is that:
$a^2_i=b^2_j=0$ for all $i, j$, that $a^1_i=0$ for all $i$ as well, and that  $b^1_j =h_i\al_l-\al_ih_l$ for some $h_i$ and some $h_l$ when $[Y_l, Y_i]=Z_j$.

Therefore, the number of free parameters for $\betal\in \Ker\delta^1$ is $q$ for $\betal(X_1)$ and $p+1$ for $\betal(X_2)$. %If $\frak N$ is free, then $p'=p+1$ And for $\betal\in \rm{Im}\delta^0$ there are $p$ free parameters for $\betal(X_1)$ and no free parameters for $\betal(X_2)$.  
So the dimension of the constant cohomology $H_\rho^1(\frak N)$ is $p+q+1$ and the cocycle family which parametrizes $H_\rho^1(\frak N)$ is defined by:

$$\betal_{a, b}(X_1)=\sum_{i=1}^q a^1_iY_i$$ 

$$ \betal_{a, b}(X_2)= \sum_{i=1}^{q} \mu \alpha_iY_{i} +\sum _{j=1}^p b^2_j Z_j$$

The part $ \sum_{i=1}^{q} \mu \al_iY_{i} $ reflects $\rho_{\al, \beta}(H_{id}^1(\frak g, \frak g))$ (where $\frak g$ is $\mathbb R^2$), namely it reflects the coordinate changes of the action $\rho_{\al, \beta}$ given by the matrix 
\begin{equation}\label{coordchange}
\pi_\mu:={\small \left(
\begin{array}{cc}
 1 & 0\\
 \mu& 0\\
\end{array}
\right)}
\end{equation} acting on $(X_1, X_2)$.

The constant family of actions $\rho_{a, b}$ for which Theorem \ref{main} gives local transversal rigidity,  is given by the cocycle family $\betal_{a, b}$, modulo the coordinate changes parametrized by the one-dimensional parameter $\mu$. Namely, $\{\rho_{a, b}\}$  is the $q+p$-dimensional family and each $\rho_{a, b}$ is defined by  commuting pair of vector fields $\sum_{i=1}^q a^1_iY_i$ and $\sum _{j=1}^p b^2_j Z_j$ where $a_i^1$ and $b_j^2$ are arbitrary constants.

\subsection{Cohomological conditions for transversal local rigidity}\label{general-result}

In \cite{D}, we obtained a general result which, given that certain cohomological conditions over a given action are satisfied,  implies transversal local rigidity. We repeat bellow the statement of that result. We shell prove that the commuting actions on 2-step nilmalifold considered in previous sections do satisfy these cohomological conditions. 

\begin{thm}\label{IFT}
Let $M$ be a smooth compact manifold, let $V=Vect^\infty(M)$ and let $\frak g$ be the Lie algebra of a Lie group $G$.

Let $\{\rho_{\mu, \la}\}_{\mu\in \RR^q, \la\in \RR^k}\in LieHom(\g, V)$ be a $C^1$-family of smooth actions on $M$ such that there exists a positive constant $\si$ and a collection $\mathcal C=\{C_r\}_r$ of positive constants for which the following conditions are satisfied for all $\mu$ in some neighborhood $\mathcal V$ of 0 in $\mathbb R^q$.

(i) The coboundary operator $\delta^0_{\rho_{\mu, 0}}$ has a $(\mathcal C, \si)$-tame inverse $\delta^{0*}_{\rho_{\mu, 0}}$ defined on the  image of $\delta^0_{\rho_{\mu, 0}}$. This means:  for every $\betal\in Im (\delta_{\rho_{\mu, 0}}^0)$ there exists $H\in \Vi $ such that $\delta_{{\rho_{\mu, 0}}}^0 H=\beta$ and $\|H\|_r\le C_r\|\beta\|_{r+\sigma}$ for all $r$.

(ii) The second coboundary operator $\delta^1_{\rho_{\mu, 0}}$ has a $(\mathcal C, \si)$-tame inverse $\delta^{1*}_{\rho_{\mu, 0}}$ on its image. 
Denote by $\Delta^\mu$ the operator $I-\delta^{1*}_{\rho_{\mu, 0}}\delta^1_{\rho_{\mu, 0}}$.

(iii) Denote by $p^{\mu}$ the projection map from  $Z^1_{\rho_{\mu, 0}}(\g, V)$ to $H^1_{\rho_{\mu, 0}}(\g, V)$ and let $P^\mu:=p^\mu\Delta^\mu$. Let $s^\mu(\mu_1, \la):=\rho_{\mu+\mu_1, \la}-\rho_{\mu,0}$. Let $P_1^{\mu}=\pi_1\circ P^{\mu}$ and $P_2^{\mu}=\pi_2\circ P^{\mu}$ denote compositions with coordinate maps, where $\pi_1$ projects $\RR^{q+k}$ to the first $q$ coordinates in $\RR^{q}$ and $\pi_2$ projects $\RR^{q+k}$ to the last $k$ coordinates, in $\RR^{k}$.  

%The main assumption on the family $\{\rho_{\mu, \la}\}$ is that it parametrizes $H^1_{\rho_{\mu, 0}}(\g, V)$, more precisely, assume that:  

 \hspace{0.2in} a) $H^1_{\rho_{\mu, 0}}(\g, V)\simeq\RR^{q+k}$ and  $s^\mu$ is a section map for $P^\mu$, namely $P^\mu\circ s^\mu=Id_{\RR^{q+k}}$. 
 
 \hspace{0.2in}  b) The maps $p^\mu$ and $s^\mu$ are  bounded, uniformly in $\mu\in \mathcal V$ with respect to the $C^0$ norm on $\frak C^1(\g, V)$ and the usual norm on $\RR^{q+k}$.
 Moreover $\|s^\mu(\mu_1, 0)\|_{\si+1}\le C\|\mu_1\|$ for all $\mu\in \mathcal V$ and the map $\mu\mapsto P_2^{\mu}$ is bounded on $\mathcal V$, with respect to the operator norm. 
 
  \hspace{0.2in}   c) Smoothing operators can be chosen so that each $P^{\mu}$ is equivariant under the action of smoothing operators (see Section 2.2 of \cite{D} for definition of smoothing operators).
  
\vspace{0.1in}
  
Then $\{\rho_{0, 0}\}$ is transversally locally rigid with respect to the family $\{\rho_{0, \la}\}$. 

%namely:  there exists $\eps>0$, $l>0$ and $B\subset \RR^{k}$, such that for every perturbation $\{\tilde\rho_{0,\la}\}$ of the family $ \{\rho_{0,\la}\}$ which is $\eps$-close to $ \{\rho_{0,\la}\}$ in $C^{l, 1, B}$-distance, there exist parameters $\bar \mu$ and $\bar \la$ such that $\tilde\rho_{0,\bar \la}$ is in the conjugacy class of $\rho_{\bar \mu,0}$ via a smooth conjugation which is $\eps$-close to the identity. 

\end{thm}

\subsection{Proof of Theorem \ref{main}}\label{proof-main}
The acting Lie group $G$ in our case is $\RR^2$, with $\frak g=Lie(G)$,  the manifold $M$ is $\GN$, and the actions $\rho_{\mu, \la}$ are given by Lie algebra homomorphisms defined on the basis of Lie algebra $\frak g$ by $$\frac{\partial}{\partial x_1}\mapsto \sum_i (\al_i+\la_i^1)Y_i$$

$$\frac{\partial}{\partial x_2}\mapsto \sum_i \mu \al_iY_i+\sum_j (\beta_j+\la_j^2)Z_j$$
The parameter $\mu$ is one-dimensional and the parameter \\
$\lambda=(\lambda_1^1, \dots, \lambda_q^1, \lambda _1^2, \dots, \lambda_p^2)$ is $q+p$-dimensional. For simplicity we may also write $\lambda=(\lambda_i^1, \lambda _j^2)$

The parameter $\mu$ reflects  coordinate changes $\pi_\mu$ of the $\mathbb R^2$ action $\rho_{0, 0}$, defined by $ \frac{\partial}{\partial x_1}\mapsto X_1:=\sum_i\al_i Y_i$ and $\frac{\partial}{\partial x_2}\mapsto X_2:= \sum_j\beta_j Z_j$, where $\pi_\mu$ are given by $\pi_\mu :X_1\mapsto X_1$ and $\pi_\mu: X_2\mapsto X_2+ \mu X_1$.

 From Sections \ref{1coh-inverse} and Section \ref{2coh-inverse} we have that both first and second coboundary operators over $\rho_{0, 0}$ with coefficients in $\CGN$ have tame inverses and that the first cohomology with coefficients in smooth functions  is trivial. By Section 5.3, Proposition 1 in \cite{D}, it follows that the same is true for coboundary operators over $\rho_{0, 0}$ with coefficients in vector fields, and that $H^1_{\rho_{0,0}}(\frak g, \rm Vect^\infty(M))=H^1_{\rho_{0,0}}(\frak g, \frak N)$. This is a consequence of the fact that $\GN$ is parallelizable and that with respect to a basis of $\frak N$ the coboundary operators have upper triangular form. 
 
 Notice that for any action $\rho$ given by a Lie algebra homomorphism from  $\frak g$ into $\frak N$,  if $\pi\in LieAut (\frak g)$, then $H^1_{\rho\circ \pi}(\frak g, V)=H^1_\rho(\frak g, V)\circ \pi^{-1}$. Since $\rho_{\mu, 0}=\rho_{0,0}\circ \pi_\mu$ (see \eqref{coordchange}), as long as for $\mu$ in a fixed neighborhood of 0, $\pi_\mu$ has a norm bounded by a constant (which is obviously true here), it follows immediately from  the corresponding conclusion for $\rho_{0,0}$, that both first and second coboundary operators over $\rho_{\mu, 0}$ with coefficients in vector fields have tame inverses. Thus  conditions (i) and (ii)  of Theorem \ref{IFT} are satisfied.

 Moreover, the same remark that  $H^1_{\rho\circ \pi}(\frak g, \frak N)=H^1_\rho(\frak g, \frak N)\circ \pi^{-1}$ also implies $\rm dim H^1_{\rho_{0,0}}(\frak g, \frak N)=\rm dim H^1_{\rho_{\mu,0}}(\frak g, \frak N)$ and consequently  
 the dimensions of all the spaces $H^1_{\rho_{\mu,0}}(\frak g, \rm Vect^\infty(M))$ are the same and they coincide with the dimension of $H^1_{\rho_{0,0}}(\frak g, \frak N)$ which is computed in Section \ref{constant-coh}. This shows the first part of condition (iii) a).

 It remains to check  the other conditions in  (iii).   
 %Part a) of condition (iii) is proved in Section \ref{constant-coh} for $\rho_{0,0}$, and follows for $\rho_{\mu, 0}$ by the same remark used above, i.e. the fact that  $H^1_{\rho\circ \pi}(\frak g, \frak N)=H^1_\rho(\frak g, \frak N)\circ \pi^{-1}$. 
 
 From  Section \ref{2coh-inverse} we see that the operator $\Delta^0=I-\delta^{1*}_{\rho_{0,0}}\delta^{1}_{\rho_{0,0}}$ on $C^\infty$ valued cocycles  takes a map $\betal: \frak g\to \CGN$ into a coboundary plus a constant. Thus the corresponding operator $\Delta^0$ on vector fields valued cocycles takes  $\betal: \frak g\to \rm Vect^\infty(\GN)$ 
 %(OR IT SHOULD BE $\beta: \rho_{0, 0}(\frak g)\to \rm Vect^\infty(\GN)$ INSTEAD???) 
 into a coboundary over $\rho_{0, 0}$ plus a constant cocycle which takes values in constant vector fields. Similar fact holds for  $P^\mu$. From Section \ref{constant-coh} the form of  these constant cocycles is known,  so for the map $P^\mu$ we have $P^\mu(\betal)(\frac{\partial}{\partial x_1})=(a_1^1, \dots, a_q^1)$ where $a_i^1$ is average of the component of $\betal(\frac{\partial}{\partial x_1})$ in the direction of $Y_i$; and $P^\mu(\betal)(\frac{\partial}{\partial x_2})= (\mu_1, b_1^2, \dots, b_p^2)$, where $\mu_1$ is the average of the component of $\betal(\frac{\partial}{\partial x_2})$  in the direction of $Y_i$ divided by $\al_1$, and $b_j$ is the average of the component of $\betal(\frac{\partial}{\partial x_2})$  in the direction of $Z_j$.  This implies the identity in part a) of condition (iii). 
 
 %Instead of $\frac{\partial}{\partial x_i}$ can write $\pi_\mu(X_i)$...............). 

 The boundedness requirement in part b) in this situation holds trivially with respect to the $C^0$ norm because the operator $p^\mu$ amounts to taking averages and $s^\mu$ takes values in constant ($\frak N$-valued) cocycles. 
 
 The map $s^\mu(\mu_1, 0)$ takes $\frac{\partial}{\partial x_1}\mapsto 0$ and $ \frac{\partial}{\partial x_2}\mapsto \mu X_2$, so the estimate required in part b) is immediate.
 
 For the part c) of (iii) we only need to know that one can choose smoothing operators which do not affect averages along various directions, and this is a well known fact, for example one possible choice of smoothing operators is described in \cite{Zehnder}, or in \cite{Hamilton}, or one can have a construction more specific to the manifold like in \cite{DK-KAM}.

 %Finally part c) is trivially satisfied since $P^\mu$ is trivial. 
 
 Thus the conclusion of the Theorem \ref{IFT} follows, namely the action $\rho_{0, 0}$ (denoted in Theorem \ref{main} by  $\rho_{\alpha, \beta}$ with $\alpha$ and  $\beta$ Diophantine) is transversally  locally rigid with respect to the family $\{\rho_{0, \la}\}$ (which in the initial statement of the Theorem \ref{main} was labeled by $\rho_{a, b}$).

\section{Further remarks on actions with GH leafwise Laplacian}\label{further_remarks}

\subsection{Higher step nilmanifolds} Let $N$ be nilpotent Lie group of step $r$ and let $\frak N$ be its Lie algebra. Let $\frak N_j=[\frak N, \frak N_{j-1}]$ , $j=1, \dots, r$ denote the lower central series of $\frak N$.  Let $\GN$ be a compact nilmanifold. 

A system of constant coefficient vector fields $\{X_1,  \dots, X_k\}$ is called globally hypoelliptic (GH) if whenever the system of equations $X_1u=f_1, \dots, X_ku=f_k$ has a distributional solution $u$ for $C^\infty$ functions $f_1, \dots, f_k$, then $u$ is also $C^\infty$. 

Let $\rho$ denote the action generated by $X_1, \dots, X_k$ and let $L_\rho$ denote the leafwise Laplacian $\sum_{i=1}^k X_i^2$. It is easy to see that if $L_\rho$  is GH then the system $\{X_1, \dots, X_k\}$ is GH. 

In \cite{CR} and \cite{CR1} Cygan and Richardson conjecture the following: 
\begin{conj} \label{conjCR} The system  $\{X_1,  \dots, X_k\}$ is GH if and only if the following two conditions hold:

i) The system is GH on the associated torus.

ii) Let $\frak  L$ denote the Lie subalgebra spanned by $X_1, \dots, X_k$. For each non-zero integral functional $\la\in (\frak N_j/\frak N_{j+1})^*$, 
\begin{equation}\label{higher-step}
\la(\frak L\cap \frak N_j+\frak N_{j+1})\ne 0, \, \, \, \, j=1, \dots, r-1
\end{equation}
($\la\in \frak N_j^*$ is integral if $\la(\log \Gamma\cap \frak N_j)\subset \mathbb Z$. ) 
\end{conj}

This conjecture is proved in \cite{CR} under additional condition that for every infinite dimensional representation $\pi$ in $\reps_\infty$ the corresponding coadjoint orbit is either flat  or $\pi$ is inducable from a polarization of codimension one in $\frak N$. In particular, the conjecture holds for 2-step nilmanifolds, since all orbits are flat, and also holds for any nilmanifolds of higher step which have all orbits flat.

Whenever the Conjecture \ref{conjCR} holds, it follows that if an $\mathbb R^k$ homogeneous action on $\GN$ has GH leafwise Laplacian then the rank of the action $k$ must be at least the number of steps $r$ of $N$. Otherwise from the condition $ii)$ above, the system $\{X_1, \dots, X_k\}$  is not GH, thus the corresponding Laplacian is not GH. 

So it seems that it may be possible to construct examples of $\mathbb R^r$ actions with GH Laplacian on $\GN$ where $N$ is of step $r$, where condition $ii)$ would be enhanced with additional Diophantine type condition on each $\frak L\cap \frak N_j$ in order to ensure that $\Ker L_{\rho}$  is GH and that tame estimates can be obtained for the inverses of cohomological operators. Moreover, such examples would have trivial first and second cohomology. Namely in \cite{CR} it is also proved that condition $ii)$ also implies that $\Ker L_{\rho}$ for action $\rho$ spanned by $X_1, \dots, X_k$, contains only constants i.e. is trivial. So it is possible that the transversal local rigidity result may be extended to such examples.

%\subsection{Globally hypoelliptic systems of vector fields}

\subsection{On the existence of $\R^k$ algebraic actions with globally hypoelliptic leafwise Laplacian} 
The condition on global hypoellipticity of the leafwise Laplacian $L_\rho$ of an abelian action $\rho$ is obviously a very strong condition, although probably not as strong as the condition on global hypoellipticity of a single vector field. In the spirit of the Conjecture \ref{GWK}, the following is proposed by Federico Rodriguez-Hertz:
\begin{conj}\label{conj-hertz}
Let  $\rho$ be  an $\R^k$ homogeneous action on the homogenenous space $\Gamma\setminus G$, where $G$ is a Lie group and $\Gamma$ a cocompact lattice in $G$. If $\rho$ has a globally hypoelliptic leafwise Laplacian $L_\rho$, then $\Gamma\setminus G$ is a nilmanifold or an infranilmanifold. 
\end{conj}
Work  towards proving this conjecture is in progress, by the author and  J. Tanis.

\subsection{ Parametric local rigidity}

Let $\{\rho_\la\}_{\la\in \mathbb R^d}$ be a finite dimensional $C^1$-family of smooth $\mathbb R^k$ actions on a smooth manifold $M$. We say that $\rho_0$ is \emph{parameter locally rigid} with respect to  the family $\{\rho_\la\}$ (or locally rigid modulo the family $\{\rho_\la\}$), if for every sufficiently small perturbation $\tilde \rho_0$ of $\rho_0$ there exists a parameter $\bar \la$ close to $0$ such that $\tilde \rho_0$ is in the smooth conjugacy class of $\rho_{\bar\la}$. 

The classical result Theorem \ref{KAMtorus} states that every flow generated by a Diophantine vector field on the torus is parameter locally rigid with respect to the family of actions generated by the constant vector fields on the torus. 

In the context of the Theorem \ref{main} on a 2-step nilmanifold (which can be viewed as being fibered over a torus), this raises a question whether it is possible to strengthen the transversal local rigidity result in Theorem \ref{main} to claim parametric local rigidity. Possibly under additional assumptions it may be that the Diophantine condition on the base torus and the Diophantine condition on the torus in the 0-fiber of the nilmanifold, may serve as moduli of smooth conjugacy for \emph{arbitrary} smooth perturbations of $\rho_{\al, \beta}$.
%Perhaps with the inv. measure assumption--rotation numbers on the torus..--> Bassam's agrument

\end{document}